\documentclass[english, reqno]{amsart}
\usepackage{babel}

\usepackage{setspace}

\usepackage{mathtools}

\usepackage{amsmath ,amsthm,amssymb,babel,amsfonts,graphics}

\newtheorem{thm}{Theorem}[]

\newcommand{\del}{\mbox{$\partial$}}

\newcommand{\half}{\mbox{$\frac{1}{2}$}}

\newcommand{\IR}{\text{${\mathbb{R}}$}}

\newcommand{\IZ}{\text{${\mathbb{Z}}$}}

\newcommand{\IN}{\text{${\mathbb{N}}$}}

\newcommand{\p}{\mathfrak{p}}
\newcommand{\q}{\mathfrak{q}}
\newcommand{\s}{\mathfrak{s}}

\newcommand{\comment}[1]{}

\title{On the Hartree--Fock Eigenvalue Problem}

\author{R. A. Zalik}

\address{Department of Mathematics and Statistics, Auburn University,
Auburn, Al. 36849--5310}\email{zalikri@auburn.edu}

\thanks{The author is grateful to Professor Vincent Ortiz, Department of Chemistry and Biochemistry, Auburn University, for his helpful comments.}
 \thanks{arXiv2408.07708}

\date{}

\begin{document}
\doublespacing

\begin{abstract}
A differential equation may be transformed into an algebraic equation by convolution with another function. In this paper we focus on this approach.
Since there are a number of functions that may be suitable, the convolution approach may lead to alternative ways of finding approximate solutions of the Roothaan--Hall equations.

 \end{abstract}

\maketitle

\pagestyle{myheadings}
\markboth{ Hartree--Fock}{ZALIK}

\section{Introduction} \label{section1}

The  Hartree--Fock  equations are nonlinear eigenvalue equations of the form
\begin{multline} \label{eq1}
\left [-\half\nabla^2 - \sum_{j=1}^n \frac{Z_j}{|x-\xi_j|} \right] \psi_a(x)+2\sum_{c=1}^n \int |\psi_c(s)|^2 \frac{1}{|x-s|} \,ds\,\psi_a(x)\\
-\sum_{c=1}^n \int \psi_c^{\ast}(s)\psi_a(s) \frac{1}{|x-s|}\, ds\, \psi_c(x)=\varepsilon_a \psi_a(x),\\
x\in\IR^3, \quad a=1, \dots, n,
\end{multline}
where $Z_j$ is the atomic number of nucleus $N_j$, which is located at point $\xi_j \in \IR^3$, and $\nabla^2$  is the Laplace operator.
These equations apply only to atoms or molecules that are closed shell systems with all orbitals doubly occupied. Thus they have an
even number $N$ of electrons, and $n=N/2$; they are the basic approximate solutions of the electronic Schr\"odinger equations. These solutions are functions $\psi_a$, called spin orbitals, with orbital energies $\varepsilon_a$;  they are continuous on $\IR^3$ and are usually assumed to form an orthonormal set in $L_2$. Since $|\psi_a|^2$ is a probability measure we also assume that  $|\psi_a(x)| \le 1$. The $N$ spin orbitals with lowest energy are called \emph{occupied} spin orbitals. 

Equations \eqref{eq1} imply that the functions $\nabla^2\psi_a$ are defined everywhere on $\IR^3$ except perhaps at the points $\xi_j$. Since the kinetic energy of the electrons must be finite, we will also assume that the partial derivatives $\del^2\psi_a/\del x_j^2, \, j=1,2,3$,  are essentially bounded.  The Gaussian--type functions used in Quantum Chemistry satisfy these conditions.

This eigenvalue problem, although similar, is not of the type studied by Kato: the function $W$ defined in \cite[(2.2)]{Ka} is expressed in terms of a bounded function $W_0$ and functions $V_{ij}$ of compact support, and therefore that problem does not match the pattern of \eqref{eq1}.

The method of solution of \eqref{eq1} pioneered by Roothaan and Hall consists of expanding the orbitals by means of a suitable basis, transforming the Hartree--Fock equations into a nonlinear system of algebraic equations, and then trying to solve this system. These and related topics are studied in standard introductory texts in Quantum Chemistry such as \cite{IL,SO}, advanced texts such the classical \cite{MS, Ra}, more recent books such as L. Lin and J. Lu \cite{LL2}, specialized works such as \cite{Da},
and in many research articles.  It is also discussed in a number of classical treatises on Quantum Mechanics such as von Neumann \cite{vN}, Gottlieb and Yan \cite{GY}, and Landau and Lifschitz \cite{LL}. Their numerical solution is implemented in several software packages.
The Roothaan--Hall  method implicitly assumes that the basis expansion can be differentiated term by term; another consideration that should influence the choice of a basis follows from the virial theorem, which states that the expectation value of the potential energy operator equals two times the expectation value of the kinetic energy operator. This theorem is satisfied by exact wave functions as well as by the bases used in practice, which contain exponential functions (but it would not be satisfied, for example, by bases containing Haar wavelets). The systems of nonlinear equations created using these bases  are solved using an ansatsz which is computationally expensive and does not always converge; so it is worthwhile exploring alternative approaches. 
The alternative we shall explore is the use of convolutions, which will allow us to transform \eqref{eq1} into a purely algebraic problem, in which the Laplacian is eliminated. This approach does not seem to have been previously studied. Our motivation can be found in the books by Folland \cite{Fo, Fo1}, where the author uses convolution with the Poisson kernel for a different but similar purpose, namely to transform partial differential equations into ordinary differential equations (see Theorem \ref{thm5} below).

Our point of view is strictly mathematical: apart from the initial hypotheses outlined in this preface, we shall have no recourse to physical intuition or experimental evidence. 

\section{Notation} \label{section 2}
We will use standard mathematical notation. In particular, 
$\IZ$ will denote the set of integers, $\IZ^+$ the set of strictly positive integers, $\IN$ the set of nonnegative integers, and $\IR$ the set of real numbers. 
Unless otherwise indicated $k \in \IZ; m,n \in \IZ^+; \eta, \xi \in \IR;  x=(x_1,x_2,x_3), s=(s_1,s_2,s_3) \in \IR^3;
||x||=(x_1^2+x_2^2+x_3^2)^{1/2}$;  if $z$ is a complex number, $z^{\ast}$ will denote its complex conjugate. We set
\[ 
h(x)=|x|^{-1}, \quad h_c(x)=h(x-c), \quad \int = \int_{\IR^3},
\]
\[
\quad\del_t=\frac{\del}{\del t}, \quad \del^2_t=\frac{\del^2}{\del t^2},  \quad \del^n_j= \left (\frac{\del}{\del_{x_j}} \right )^n,
\]
and if $\alpha=(\alpha, \dots, \alpha_n)\in \IN^n$, then 
\[
 \del^{\alpha}= \left(\frac{\del}{\del x_1} \right)^{\alpha_1} \cdots \left(\frac{\del}{\del x_n} \right)^{\alpha_n}.
\]

We say that $f$ is in $L_{\infty}(\IR^3)$ or $L_{\infty}$ if there is a constant $K$ such that
$|f(x)| \le K$ for almost all $x \in \IR^3$. The $L_{\infty}$ norm is defined as
\[
||f||_{\infty} = \text{ess sup}|f(x)|, \qquad x \in \IR^3.
\]
$C_0$ is the set of continuous functions on $\IR^3$ that vanish at infinity; for $f \in C_0$
\[
||f||_u=\max |f(x)|,  \qquad x \in \IR^3.
\]
$C_{\infty}$ is the set of infinitely differentiable functions in $\IR^3$.

If $S \subset  \IR^3$ and $1 \le p <\infty$ we say that $f$ is in $L_p(S)$ if
\[
\int_S |f(x)|^p\,dx < \infty,
\]
and  $L_p(\IR^3)$ will be abbreviated $L_p$. The
$L_p(\IR^3)$ norm of $f$ is defined by
\[
||f\||_p =\left (\int |f(x)|^p \, dx\right )^{1/p}.
\]
If $f,g \in L^2$,
\[
\langle f,g \rangle = \int f^{\ast}(x) g(x)\, drx
\]
The convolution $f\ast g$, of $f$ and $g$, is defined by
\[
(f\ast g)(x)=\int f(s) g(x-s)\, ds =\int f(x-s) g(s)\, ds.
\]

$\psi_a$ will always denote a spin orbital.

Additional notation will be introduced as needed.
\section{Preliminary Results} \label{section 3}
\begin{thm} \label{thm1} Let $V$ denote the set of points in $\IR^3$ other than the singular points $\xi_c$; then 
$\nabla^2\psi_a \in C(V) \cap L_2$. If $\psi_a$ is also in $L_1$ (or if $\lim_{x\rightarrow \infty}\psi_a(x)=0$), then $\nabla^2\psi_a $ is in $L_1$ as well.
\end{thm}
\begin{proof}
In \cite{Za1} we showed that if $\psi_a \in L_1$ then $\lim_{x\rightarrow \infty} \psi_a(x)=0$.

For $a=1, \dots, n$, let 
\begin{equation*} 
\begin{split}
\s_{a c}(x)= &\int \psi_c^{\ast}(s)\psi_a(s) h(x-s)\,ds=[(\psi_c^{\ast} \psi_a)\ast h](x),\\
\p(x)= &2\sum_{c=1}^n Z_c h(x-\xi_c)=2\sum_{c=1}^n Z_c h_{\xi_c}(x),  \quad\text{and}\\
\q(x)=  & 4\sum_{c=1}^n \int |\psi_c(s)|^2h_{\xi_c}(x)\,ds  = 4\sum_{c=1}^n  s_{c,c}(x).\\
\end{split}
\end{equation*}
Since \eqref{eq1} is equivalent to
\begin{equation} \label{eq2}
\nabla^2\psi_a(x)=[\q(x)-\p(x) -\varepsilon_a]\psi_a(x)-2\sum_{c=1}^n\s_{a,c}(x)\psi_c(x)
\end{equation}
we will focus our attention on the functions  $\p\psi_a$, $\q\psi_a, \psi_a$, and $\s_{a,c}\psi_c$. We know that
\[
\psi_a \in L_{\infty}\cap L_2 \cap C(\IR^2),
\]
that the definition of $\p$ implies that it is in $C(V)$, that $\q$ is expressed in terms of the functions $\s_{c,c}$, and that
$\p$ is expressed in terms of the functions $h_{\xi_c}$.  We will now show that
\begin{equation} \label{eq3}
\s_{a,c} \in L_{\infty} \cap C(V), \quad  h_{\xi_c} \psi_a \in L_2 ,
\end{equation}
and, if $\psi_a \in L_1$, then $h_{\xi_c} \psi_a \in L_1$ as well.

We begin with $h_{\xi_c} \psi_a$. Let $\tau \in \IR^3$; then
\begin{multline} \label{eq4}
||h_{\tau}\psi_a||^2_2=\int |h_{\tau}(s)\psi_a(s)|^2 \,ds= \\
\int |s-\tau|^{-2}\,|\psi_a(s)|^2 \,ds=\int |s|^{-2}\, |\psi_a(s+\tau)|^2 \,ds =\\
\int _{|s|\le1}|s|^{-2}|\psi_a(s+\tau)|^2 \,ds + \int_{|s| >1} |s|^{-2}|\psi_a(s+\tau)|^2 \,ds\le\\
\left (\int _{|s|\le1}|s|^{-2}\,ds \right)^{\half} \left (\int _{|s|\le1}|\psi_a(s+\tau)|^2 \,ds \right )^{\half} +\int_{|s| >1}|\psi_a(s+\tau)|^2 \,ds \le \\
4\pi+1.
\end{multline}
On the other hand, if $\psi_a \in L_1$,
\begin{multline*}
||h_{\tau}\psi_a||_1= \int |h_{\tau}(s)\psi_a(s)|\,ds= \int |s|^{-1}\, |\psi_a(s+\tau)| \,ds =\\
\int _{|s|\le1}|s|^{-1}|\psi_a(s+\tau)| \,ds+\int _{|s| > 1}|s|^{-1}|\psi_a(s+\tau)| \,ds \le\\
\left ( \int _{|s|\le1}|s|^{-2} \,ds\right )^{\half} \left ( \int _{|s|\le1}|\psi_a(s+\tau)|^2 \,ds \right )^{\half} + \\
\left ( \int _{|s| >1}|s|^{-2} \,ds\right )^{\half} \left ( \int _{|s|>1}|\psi_a(s+\tau)|^2 \,ds \right )^{\half} \le\\
\left ( \int _{|s|\le1}|s|^{-2} \,ds\right )^{\half}+\left ( \int _{|s|>1}|\psi_a(s+\tau)|^2 \,ds \right )^{\half} =
2\sqrt{\pi} +1.
\end{multline*} 
We now turn to $\s_{a,c}$.
\begin{multline*}
|\s_{a,c}| \le
 \int |\psi_a(s) \psi_c(s)h(x-s)|\,ds = \int |\psi_a(x-s) \psi_c(x-s)| \, |s|^{-1}\,ds=\\
 \int_{|s|\le 1}|\psi_a(x-s) \psi_c(x-s)| \,|s|^{-1}\,ds + \int_{|s| > 1}|\psi_a(x-s) \psi_c(x-s)| \,|s|^{-1}\,ds=
 I_1(x)+I_2(x).
\end{multline*} 
Applying H\"older's inequality and switching to spherical coordinates we get
\begin{multline*}
I_1(x) \le \left ( \int_{|s|\le 1}|\psi^2_a(x-s) \psi^2_c(x-s)|\,ds \right )^{\half} \left (\int_{|s|\le1}|s|^{-2}\, ds \right)^{\half} \le \\
\left (\int_{|s|\le1}|s|^{-2}\, ds \right)^{\half} =2\sqrt{\pi},
\end{multline*}
and we readily conclude that
\begin{equation} \label{eq5}
|\s_{a,c}| \le 2\sqrt{\pi}+1.
\end{equation}
To show that $\s_{a,c}$ is continuous on $V$ we use a similar procedure; for convenience we set $g(x)=\psi_a(x)\psi_c^{\ast}(x)$. 
Let $\varepsilon > 0$ and $M=4/\varepsilon$; then
\begin{multline*}
|\s_{a,c}(x+\tau)-\s_{a,c}(x)| \le \int |g(x+\tau-s)-g(x-s))| \, |s|^{-1} \, ds= \\
 \int_{|s|\le M} |g(x+\tau-s)-g(x-s))| \, |s|^{-1} \,ds + \int _{|s| > M}|g(x+\tau-s)-g(x-s))| \, |s|^{-1} \,ds =\\
 I_3(x)+I_4(x).
\end{multline*}
Let $K= (4 \pi M^2)^{-1}$. Since the spin orbitals are continuous on 
$\IR^3$ it follows that $g(x)$ is continuous on $\IR^3$ and therefore uniformly continuous on compact sets. Since the set of singular 
points $\xi_c$ is finite, $V$ is an open set; let $x \in V$, then there is a $\delta > 0$ such that if $|\tau| < \delta$ then $x+\tau \in V$ and
\[
|g(x+\tau-s)-g(x-s)| <K\varepsilon.
\]
Hence
\[
I_3 \le K\varepsilon \int_{|s| \le M}|s|^{-1}\, ds=K\varepsilon (2\pi M^2)=\frac{\varepsilon}{2}.
\]
We also have
\begin{multline*}
I_4(x) \le \frac{\varepsilon}{4} \int |g(x+\tau-s)+g(x-s))| \,ds \le \\
\frac{\varepsilon}{2} ||g||_1= \frac{\varepsilon}{2}||\psi_a \psi_c||_1 \le \frac{\varepsilon}{2}||\psi_a||_2 ||\psi_c||_2=\frac{\varepsilon}{2},
\end{multline*}
Therefore 
\[
|\s_{a,c}(x+\tau)-\s_{a,c}(x)| \le \varepsilon,
\]
and the continuity of $\s_{a,c}$ on $V$ follows.

\end{proof}
Inequality \eqref{eq5} was proved by the author in \cite[Lemma 1]{Za1}; it is included here for the sake of completeness.

To transform the Hartree--Fock equations into a (nonlinear) algebraic system of equations we need to commute these operators, which reduces to differentiating functions under the integral sign; but since $\psi_a$ may not be in $L_1$ we cannot use standard theorems to do so. To motivate the necessity of such theorems consider the following elementary example: let $t \in \IR$, let $u(t)=t$ for $0 \le t \le 1$ and $u(t)=0$ otherwise, and let $v(t)=t$; then $(u \ast v)(\eta)=\eta/2 -1/3$, whence $(u \ast v)^{\prime}(\eta)=1/2$. Similarly, $(u \ast v^{\prime})(\eta)=1/2$, as well, but on the other side $(u ^{\prime}\ast v)(\eta)= \eta-1/2$. Note that $v$ is not in $L_1$. More complicated examples could be constructed that more closely match our situation, but it is probably not worth the effort.

To prove the next result we need to introduce additional notation.

Let  $e_1=(1,0,0), e_2=(0,1,0), e_3=(0,0,1)$, and let $f(x)$ be any function in $\IR^3$. We define $f_j(x;\eta)=f(x +\eta e_j)$, and $f(x;\eta)=\sum_{j=1}^3 f_j(x;\eta)$. The expression $f_j[x;0, \eta, 2\eta]$ will denote the divided difference of order 2 of $f_j (x;\cdot)$ at $0, \eta$, and $2\eta$; thus
\[
f_{j}[x;0, \eta, 2\eta]=\frac{f(x)-2f(x+\eta e_j) + f(x+2\eta e_j)}{2\eta^2}.
\]
The divided difference of order 2 of $f(x;\cdot)$ at $0, \eta$, and $2\eta$ is similarly defined.
In the next theorem $\psi$ will denote any function, not necessarily a spin orbital.
\begin{thm} \label{thm2}
Let $\psi$ and $w$ be such that their partial derivatives 
\begin{equation} \label{eq6}
\del_j\psi(x), \quad \del_jw(x), \quad \del^2_j\psi(x),\; \text{and}\quad\del^2_j w(x), \quad j=1, 2, 3,
\end{equation}
exist everywhere except, perhaps, on a set $G$ of isolated points, and assume that $\psi, w$, and the partial derivatives of order 2  displayed on \eqref{eq6} are in $L_2\cap L_{\infty}$.; assume also that $w\in L_1$. Then $[(\nabla^2\psi )\ast w](x) $and $[\psi\ast \nabla^2w](x)$ are  defined for every $x \in \IR^3$, are in $C_0$, and are uniformly continuous on $\IR^3$. Moreover\\
(a)\ If  $\del^2_j(\psi\ast w)(x)$, $j=1,2,3$,  exist everywhere except on $G$, then $\nabla^2[\psi\ast w](x)$ exists everywhere except on $G$. If, moreover,
\begin{equation} \label{eq7}
\nabla^2[\psi\ast w](x)= [\psi\ast \nabla^2w](x),
\end{equation}
then 
\[
[(\nabla^2\psi )\ast w](x)= [\psi \ast (\nabla^2 w)](x),
\]
and this identity holds for every $x \in \IR^3$.\\
(b)\ If $(1+|x)|\del_jw(x)$ and $(1+|x)|\del^2_jw(x)$ are in $L_{\infty}$ for $j=1,2,3$, then $\nabla^2[\psi\ast w](x)$ exists everywhere except on $G$, and $w(x)$ satisfies \eqref{eq7} as well.
\end{thm}
\begin{proof}
The existence of $[(\nabla^2\psi )\ast w](x)$ on $\IR^3 \setminus G$ follows directly from the hypotheses.\\
(a)\ Since the functions $\psi, w, \nabla^2\psi$, and $\nabla^2w$ are in $L_2$, we deduce that  $[(\nabla^2\psi )\ast w](x)$ and $ [\psi \ast (\nabla^2 w)](x)$ are defined for every $x \in \IR^3$ and that they are in $C_0$ and uniformly continuous (cf. \cite[Theorem 8.8]{Fo1}).

Let 
\[
M(\psi)=||\del_1^2\psi||_{\infty}+||\del_2^2\psi||_{\infty}+||\del_3^2\psi||_{\infty,}
\]
\[
D^2\psi_j(x;\eta)=\del^2_t \psi(x+t u_j) \Big{|}_{t=\eta}, \quad \text{and} \quad D^2\psi(x;\eta)= \sum_{j=1}^3 D^2\psi_j(x;\eta).
\]
If $x \notin G$, then the partial derivatives of $w$ displayed on \eqref{eq6} exist on a neighborhood of $x$; thus, from
the properties of divided differences (\cite{SB}) we know that 
\[
\lim_{\eta \rightarrow 0}\psi_{j}[x;0, \eta, 2\eta]=\half \del_j^2\psi(x) 
\]
and therefore
\begin{equation} \label{eq8}
\lim_{\eta \rightarrow 0} \psi[x;0, \eta, 2 \eta]=\half \nabla^2\psi(x).
\end{equation}
On the other hand, for $\eta$ sufficiently small there is a number $\sigma_j$ between 0 and $\eta$ such that
\[
\psi_j[x;0, \eta, 2 \eta] = D^2\psi_j(x,\sigma_j).
\]
But
\[
D^2\psi_j(x, \sigma_j)= \left ( \frac {\del^2\psi}{\del u_j^2}\right ) \Bigg{| }_{u_j=x + \sigma_j e_j},
\]
and we conclude that
\[
||D^2\psi(x;\cdot))||_{\infty} \le M(\psi).
\]
whence also
\[
\big{|} \psi[x; 0, \eta, 2 \eta] \,\big{|}  \le M(\psi).
\]
Since
\[
(\psi\ast w)[x; 0, \eta, 2 \eta]= \int \psi[x-s; 0, \eta, 2 \eta] w(s)\,ds
\]
and $w\in L_1$, applying the Bounded Convergence theorem we conclude from \eqref{eq8} that for any $x \notin G$, arbitrary but fixed,
\[
\lim_{\eta \rightarrow 0}(\psi\ast w)[x; 0, \eta, 2 \eta]= \half [\nabla^2\psi\ast w](x).
\]
But the properties of divided differences also imply that
\[
\lim_{\eta \rightarrow 0}(\psi\ast w)_{j}[x;0, \eta, 2\eta]=\half \del^2_j(\psi\ast w)(x), \quad j=1, 2 ,3,
\]
and therefore
\begin{equation} \label{eq9}
\lim_{\eta \rightarrow 0} (\psi\ast w)[x;0, \eta, 2 \eta]=\half \nabla^2(\psi\ast w)(x),
\end{equation}
whence the assertion follows.

(b)\  Let $j=1,2,3$. An application of the Bounded Convergence theorem readily shows that $\del_jw$(x) and $\del^2_jw(x)$ are defined in $x\setminus G$.Let $j=1, 2, 3$, $x$ arbitrary but fixed, and assume that the functions $\del^2_jw(x)$, $j=1,2,3$ are defined.  The hypotheses imply that there is a constant $C$ such that for $0 \le \eta \le1$. 
\[
|x|\,|D^2w_j(x, \eta)| \le C \quad a.e.
\]
and a number $\xi_j$ between 0 and $\eta$, such that
\[
w_j[x; 0, \eta,2\eta]=D^2w_j(x, \xi_j),
\]
Thus
\[
|x|\,|w_j[x; 0, \eta,2\eta]\le C _1\quad a.e.,
\]
which implies that 
\[
|x|\,|w[x; 0, \eta,2\eta]|\le 3C_1 \quad a.e.
\]
Therefore, for any $x$, $\psi(s) w[x-s; 0, \eta, 2 \eta]$ is in $L_1$.
Since
\[
(\psi\ast w)[x; 0, \eta, 2 \eta]= \int \psi(s) w[x-s; 0, \eta, 2 \eta]\,ds,
\]
applying the Bounded Convergence theorem we conclude that
\[
\lim_{\eta \rightarrow 0^+}(\psi\ast w)[x; 0, \eta, 2 \eta]= \half [\psi\ast \nabla^2]w(x),
\]
and the assertion follows from \eqref{eq9}.
\end{proof}
For any function $f$, any nonnegative integer $n$, and any multi--index $\alpha$, let
\[
||f||_{n, \alpha}= \sup_{x \in \IR^d} (1+|x|^n)|\del^{\alpha}f(x)|,
\]
and
\[
\mathcal{S}=\{r \in C_{\infty} : ||f||_{n,\alpha} < \infty \text{ for all} \; n, \alpha\}.
\]
$\mathcal{S}$ is called the \emph{Schwartz space} in $\IR^d$. In view of Theorem \ref{thm2}(b), it is clear that the functions in  Schwartz space satisfy \eqref{eq7}.
The  Gaussian--type functions used in Quantum Chemistry are in Schwartz space in $\IR^3$.
\section{Convolution} \label{section 4}
We now study the effects of convolution on the Hartree--Fock equations. 
\begin{thm} \label{thm3}
Let $w$ be such that its partial derivatives of order less or equal to 2 displayed on \eqref{eq6} exist almost everywhere, that $w$ and its partial derivatives of order 2  displayed on \eqref{eq6} are in $ L_1\cap L_2 \cap L_{\infty}$,  and that
\[
\nabla^2[\psi_a \ast w](x)= [\psi_a\ast \nabla^2w](x).
\]
Then for $x \in \IR^3$ and $a=1, \dots n$,
\begin{equation} \label{eq10}
 \varepsilon_a [\psi_a \ast w](x)=
\left [\left ( (\q - \p) \psi _a-2\sum_{c=1}^n \s_{a,c}\psi_{c}\right )\ast w\right ](x) -[\psi_{a}\ast \nabla^2w](x)
\end{equation}
and $\psi _a\ast w$ is in $C_0$  and is uniformly continuous on $\IR^3$.
\end{thm}
\begin{proof}
That the partial derivatives $\del_j\psi_a$ and  $\del^2_j\psi_a$, $j=1, 2, 3$, exist almost everywhere is implicit in \eqref{eq1}. The assumptions made in the Introduction imply that $\psi_a$ and $\del^2_j\psi_a$, $j=1, 2, 3$, are in $L_2\cap L_{\infty}$.  Convolving both sides of \eqref{eq2} with $w$ and applying Theorem \ref{thm2},  the assertion follows.
 \end{proof}
 Let $\{\phi_{\tau}; \tau \in \IN^3\}$ be a basis of $L_2$; then there are scalars $b_{\tau,a}$ such  that if
 \[
 T_{a,\nu}=\sum_{|\tau| \le \nu} b_{\tau,a}\phi_{\tau},
 \]
 then 
 \[
 \lim_{n\rightarrow \infty} ||\psi _a-  T_{a,\nu}||_2=0.
 \] 
 Since the sequence of functions $\{T_{a,\nu}; \nu \in \IN^3\}$ is convergent in $L_2$, there is a constant $K_a$ such that, for $\tau \in \IZ^3$,
$||T_{\nu,a}||_2 \le K_a$.

 We define
 \[
 q_{\nu}(x)=4\sum_{c=1}^n\s_{c,c,\nu}(x),\quad \text {where}\quad \s_{a,c,\nu}(x)=[(T^{\ast}_{a,\nu}T_{c,\nu})\ast h](x).
  \]
\begin{thm}\label{thm4}\hspace{1in}\\
(a)\ Let
\begin{equation} \label{eq11}
 \varepsilon_{a,\nu}= 
 \frac{\big{\Vert} [(\q -\p )T_{a,\nu } -2\sum_{c=1}^n(\s_{a,c,\nu} T_{c,\nu}) ] \ast w -  [ T_{a, \nu}\ast \nabla^2w ] \big{\Vert}_2}
 {||T_{\nu,a} \ast w ||_2};
 \end{equation}
then $\lim_{\nu\rightarrow \infty}\varepsilon_{a,\nu}=\varepsilon_a$.\\
(b)\ $\{[T_{\nu,a}\ast w](x); \nu \in \IR^3\}$ converges to $[\psi_a\ast w](x)$ uniformly on $\IR^3$, and
\begin{equation} \label{eq12}
|| \psi_a\ast w-T_{\nu,a}\ast w||_u \le |||\psi_a-T_{\nu,a}||_2\,||w||_2.
\end{equation}
(c)\ Assume that $[\psi_a\ast w](x_0)\ne 0$; then for $\nu$ sufficiently large $[T_{a,\nu}\ast w](x_0) \ne 0$ and, if
\begin{equation} \label{eq13}
\varepsilon_{\nu,a}(x_0)= \frac{\{[(\q _{\nu}-\p )T_{a,\nu } -2\sum_{c=1}^n(\s_{a,c,\nu} T_{c,\nu})\ast] w \}(x_0)-  [ T_{a, \nu}\ast \nabla^2w ] (x_0)}{[T_{a,\nu}\ast w](x_0)},
\end{equation}
then $\lim_{\nu \rightarrow \infty}\varepsilon_{a,\nu}(x_0)= \varepsilon_a$.
\end{thm}
\begin{proof}
(a)\ Since $\psi -  T_{\nu,a} \in L_2$, Young's inequality implies that $||T_{\nu,a} \ast w ||_2$ converges to $||\psi \ast w ||_2$ and that
the  numerator of \eqref{eq11} converges to 
\[
||[(\q\psi_{a} - \p \psi_{a} -2\sum_{c=1}^n (\s_{a,c}\psi_{c})\ast w ] -[\psi_{a}\ast \nabla^2w]||_2,
\]
and the assertion follows from Theorem \ref{thm3}.\\
(b)\ Inequality \eqref{eq12} follows from Young's inequality.\\
(c)\ That the denominator of \eqref{eq13} does not vanish for sufficiently large $\nu$, as well as the convergence of $(T_{a, \nu} \ast w)(x_0)$ to
$(\psi_a \ast w)(x_0)$, follow from \eqref{eq12}. 
Proceeding as in the proof of Theorem \ref{thm1} we see that
\begin{multline*}
||h_{\tau}T_{a,\nu} ||^2_2=\int |h_{\tau}(s)T_{a,\nu}(s)|^2 \,ds \le\\
\left (\int _{|s|\le1}|s|^{-2}\,ds \right)^{\half} \left (\int _{|s|\le1}|T_{a,\nu}(s+\tau)|^2 \,ds \right )^{\half} +\int_{|s| >1}|T_{a,\nu}(s+\tau)|^2 \,ds \le\\
4\pi K_a+K^2_a,
\end{multline*}
which implies that $||\p T_{a,\nu}||_2$ is bounded, uniformly on $\nu$, and by another application of Young's inequality
as in \eqref{eq12} above we conclude that $[(\p T_{a,\nu})\ast w](x_0)$ converges to $[(\p\psi_a)\ast w](x_0)$.
From \eqref{eq12} we know there is a constant $M_a$ such that $||T_{a,\nu}||_u \le M_a$. Repeating the arguments used in the proof of Theorem 1 we readily obtain
\begin{multline*}
|\s_{a,c,\nu}(x)| \le\\
 \left ( \int_{|s|\le 1}|T^2_{a,\nu}(x-s) T^2_{c,\nu}(x-s)|\,ds \right )^{\half} \left (\int_{|s|\le1}|s|^{-2}\, ds \right)^{\half} +\\
\int_{|s| > 1}|T_{a,\nu}(x-s) T_{c,\nu}(x-s)| \,|s|^{-1}\,ds\le\\
2 \sqrt{\pi} M_a^2+K^2_a,
\end{multline*}
which also implies the boundedness of $\{\q_{\nu}; \nu \in \IZ^3\}$ in $L_{\infty}$. The uniform bounds for 
$||h_{\tau}T_{a,\nu }||_2$ and $\s_{a,c,\nu}$ that we have just obtained are the counterparts of \eqref{eq4} and \eqref{eq5}, and using them we can readily show that the sequence 
\[
\{[( \q_{\nu}-\p )T_{a,\nu} -2\sum_{c=1}^n(\s_{a,c,\nu} T_{c,\nu}) ); \nu \in \IZ^3\}
\]
is bounded in $L_2$, and the assertion follows by convolving with $w$ and applying Young's inequality as in \eqref{eq12}.
\end{proof}

\section{The Poisson Kernel} \label{section 5}
Here is an example of a function that is not in Schwartz space, but is of sufficiently fast decay that \eqref{eq7} is satisfied.

The upper half--space in $\IR^{d+1}$ is the set 
\[
R_+^{d+1}=\{(x,t) \in \IR^{d+1}; 0 < t < \infty\},
\]
and its closure is
\[
\overline{R_+^{d+1}}=\{(x,t) \in \IR^{d+1}; 0  \le t < \infty\}.
\]
If $d \in \IN$ and $x \in \IR^d$,  the Poisson kernel for the upper half--space in $\IR^{d+1}$ is
\[
P_t(x)=\frac{\Gamma( (d+1)/2) )}{\pi^{(d+1)/2)}} {t(t^2+|x|^2)^{-(d+1)/2}}. 
\]
For $d=3$ this becomes
\begin{equation} \label{eq14}
P_t(x)=\pi^{-2} t\, (t^2+|x|^2)^{-2}.
\end{equation}
In some printings of \cite{Fo1} there is a typographical error in the definition of $P_t$.

The following is a summary of well--known properties of  the Poisson kernel.  We have culled parts (a) through (f) from \cite[Theorem 2.44]{Fo}, \cite[Theorem 2.45]{Fo}, and \cite[Theorem 8.14]{Fo1}. We include this summary for ease of reference. 
\begin{thm} \label{thm5}
Assume that $f \in L_p \; ( 1\le p \le \infty)$ and let $v(x,t)=(f\ast P_t)(x)$. Then\\
\noindent\textup(a)\ \quad$(\nabla^2 +\del_t^2)v=0 $ on $\overline{R_+^{d+1}}$. \\
\textup(b)\ \quad $\lim_{t \rightarrow 0^+} ||v(\cdot, t) - f||_p=0$, provided $1\le p < \infty$.\\
\textup(c)\ \quad $\lim_{t \rightarrow 0^+} v(x, t) = f(x)$ for a.e. $x$, and for every $x$ at which $f$ is continuous.\\
\textup(d)\ \quad If $f$ is bounded and continuous then $v$ is continuous on $\overline{R_+^{d+1}}$ and\\
 \indent\indent $v(x,0)=f(x)$.\\
\textup(e)\ \quad If $f \in L_{\infty}$ and is uniformly continuous on $R^d$, then $\lim_{t \rightarrow 0^+} v(x, t) = f(x)$ 
\indent\indent uniformly on $x$.\\
\textup(f)\ \quad If $f$ is continuous and vanishes at infinity in $\IR^d$, then \\
\indent\indent$v(x,t)$ vanishes at infinity on $\overline{R_+^{d+1}}$ and is the unique solution on $\IR_+^{d+1}$ of
\[
\nabla^2 w+\del^2_t w=0, \quad v(x,0)=f(x)
\]
\quad\quad that vanishes at infinity.\\
\textup(g)\ 
\[
\nabla^2[f\ast P_t](x)= [f \ast \nabla^2 P_t](x).
\]
\end{thm}
Part (g)  follows from Theorem \ref{thm2}(b); see also the proof of \cite[Theorem 8.53]{Fo1}.
This theorem was the motivation for our Theorem 3 and yields information about the Poisson kernel that is not implied by Theorem \ref{thm3}.
It readily yields
\begin{thm} \label{thm6}
Let  $P_t$ be defined by \eqref{eq14}, let the functions $\p, \q$, and $\s_{a,c}$ be defined as in Theorem \ref{thm1}, 
 $ u_{a,t}=\psi_a \ast P_t$, and $||u_{a,t}||_{\infty}=\sup_{x \in \IZ^3}|u_{a,t}(x)|$.
Then\\
\noindent \textup(a)\ For $a=1, \dots, n$ and $t>0$
\begin{equation*}  
\varepsilon_a (\psi_a \ast P_t)(x)=(\psi_a \ast \del_t^2P_t)(x +\{[(\q-\p)\psi_a] \ast P_t\}(x)  -\sum_{c=1}^n [(\s_{a,c} \psi_c )\ast P_t](x). \\
\end{equation*}
\textup(b)\ \indent $||u_{a,t}||_{\infty}  \le \displaystyle\frac{1}{\sqrt{2} \pi t^{3/2}}$.\\
\textup(c)\  \indent $\lim_{t \rightarrow 0^+} ||u_{a,t}-\psi_a||_2=0.$\\
\textup(d)\ \indent $\displaystyle{ \nabla^2(\psi_a \ast P_t)(x)+(\psi_a \ast \del_t^2 P_t)(x)=0;\;  t >0, u_{a,0}(x)=\psi_a(x)}$.\\
\textup(e)\ \indent If $\psi_a$  vanishes at infinity then $u_{a,t}$ vanishes at infinity in $\overline{R_+^4}$, is\\ 
\indent \indent the unique solution that vanishes at infinity of
\[
\nabla^2 w+\del^2_t w=0 \quad \text{on}\ \IR_+^{d+1}, \quad w(x,0)=\psi_a(x),
\]
\quad\quad  and $\lim_{t \rightarrow 0^+} ||u_{a,t}-\psi_a||_{\infty}=0.$\\
\textup(f)
\[
\nabla^2[\psi_a\ast P_t](x)= [f\psi_a\ast \nabla^2 P_t](x).
\]
\end{thm}
\begin{proof}  
Part (a) follows by convolving both sides of \eqref{eq2} with $\nabla^2P_t$ and applying Theorem \ref{thm5}(a).\\
Since $\psi_a \in L_{\infty}$ and $P_t \in L_1$, Part (b) follows fromYoung's inequality. \\
Part (c) follows  from Theorem \ref{thm5}(b), Part (d) follows from Theorem \ref{thm5}(f) and \eqref{eq7},
Part (e) follows from Theorem \ref{thm5}(e) and Theorem \ref{thm5}(f), and 
Part (f) follows from Theorem \ref{thm5}(g).
\end{proof}
\section{Conclusion}
Convolutions of basis functions with the Poisson kernel or other suitable functions enable elimination of the Laplacian operator in the Roothaan-Hall equations and may be possibly used in the solution of other partial differential equations. This approach also appears to be applicable to purely numerical algorithms for the Hartree--Fock equations where grids instead of basis functions are employed.

\end{document}